\documentclass{article}\begin{document}
\centerline{\bf\large Evaluations of Series Related to Jacobi Elliptic Functions\normalsize}
\vskip .4in
\centerline{N.D.Bagis}
\centerline{Stenimahou 5 Edessa, Pellas 58200}
\centerline{Edessa, Greece}
\centerline{nikosbagis@hotmail.gr}
\vskip .2in
\[
\]

\textbf{Keywords}: Hyperbolic functions; Series; Evaluations; Elliptic Functions; Special Functions;
\[
\]
\centerline{\bf Abstract}

\begin{quote}
In this article we give evaluations of certain series of hyperbolic functions using Jacobi elliptic functions theory. We also define some new functions that enable us to give characterization of not solvable class of series.

\end{quote}

\section{Introduction}

The study of infinite sums and products of certain hyperbolic functions have been a point of interest in mathematics for more than last 200 years. Works of many great mathematicians such Euler, Gauss, Jacobi, Eisenstein, Weierstrass, Abel give partial answer to these strange at first sight series. This was done by introducing and developing the theory of elliptic functions and modular forms. Later Weber, Ramanujan, Watson, Hardy, Hecke, Poincare, etc. proceeded and developed these theories further, until today, where a very large number of scientists add continuously new and very interesting results. However, some of these sums, for example one of them is 
\begin{equation}
\sum^{\infty}_{n=1}\frac{1}{e^{nx}-1}\textrm{, }x>0,
\end{equation}
where not able to evaluated with the existing theory (we call these class of sums as ''Ghost Sums''). Our main concern here is to ''evaluate'' these kind of sums by defining new functions similar to the classical elliptic functions, having not necessary double periods and overcome this problem. We also study these new functions and try to find their properties.\\        

We begin by stating and commenting some knowing results.\\
\\

Some remarkable evaluations, of closely related sums to (1) are
\begin{equation}
\sum^{\infty}_{n=1}\frac{n^{4\nu+1}}{e^{2\pi n}-1}=\frac{B_{4\nu+2}}{8\nu+4},
\end{equation}
where $\nu$ is positive integer and $B_j$ denotes the $j$th Bernoulli number.\\
Also in [11] it have been proved that
\begin{equation}
\sum_{n\geq0,n-odd}\frac{n^{4\nu+1}}{e^{n\pi}+1}=-\frac{1}{4}Q_{4\nu+1}-\frac{2^{4\nu-1}}{2\nu+1}B_{4\nu+2},
\end{equation}
where $\nu\in\{0,1,2,\ldots\}$ and $Q_{\nu}$ is defined as
\begin{equation}
Q_{\nu}=2\left(\frac{d^{\nu}}{dx^{\nu}}\frac{1}{e^x+1}\right)_{x=0}.
\end{equation}

Continuing for $\nu$ non zero integer the sums
\begin{equation}
\sum^{\infty}_{n=1}\frac{n^{-2\nu-1}}{e^{2\pi n}-1} 
\end{equation}
can evaluated from (see [2]):\\
\\
\textbf{Theorem 1.}\\
Let $a,b>0$ with $ab=\pi^2$, and let $\nu$ be any non zero integer. Then    
$$
a^{-\nu}\left\{\frac{1}{2}\zeta(2\nu+1)+\sum^{\infty}_{n=1}\frac{n^{-2\nu-1}}{e^{2an}-1}\right\}-
(-b)^{-\nu}\left\{\frac{1}{2}\zeta(2\nu+1)+\sum^{\infty}_{n=1}\frac{n^{-2\nu-1}}{e^{2bn}-1}\right\}=
$$
\begin{equation}
=-2^{2\nu}\sum^{\nu+1}_{n=0}(-1)^n\frac{B_{2n}}{(2n)!}\frac{B_{2\nu+2-2n}}{(2\nu+2-2n)!}a^{\nu+1-n}b^n,
\end{equation}
where $\zeta(s)$ is the Riemann zeta function.\\

An example of (6) for $\nu=-1$ is
\begin{equation}
\sum^{\infty}_{n=1}\frac{n}{e^{2n\pi}-1}=\frac{1}{24}-\frac{1}{8\pi}.
\end{equation}
Eisenstein and later Ramanujan consider infinite sums of the form (Eisenstein Series): 
\begin{equation}
\sum^{\infty}_{n=1}\frac{n^{2\nu-1}q^n}{1-q^n}\textrm{, }q=e^{-\pi\sqrt{r}}\textrm{, }r>0\textrm{ with  $\nu=1,2,\ldots$}.
\end{equation}
Ramanujan has given relations of how one can evaluate them using only the second and third sum of them i.e. using only the sums for $\nu=2$ and $\nu=3$, (see [2] chapters 14-15). Formulas for evaluating derivatives of Eisenstein series was also given by Ramanujan (see [2] chapter 15 Entry 13 and [3]).\\  

The function
\begin{equation}
A=\sum^{\infty}_{n=1}\frac{q^n}{n(1-q^n)}=\sum^{\infty}_{n=1}\frac{1}{n\left(e^{\pi n \sqrt{r}}-1\right)}\textrm{, }q=e^{-\pi\sqrt{r}}\textrm{, }r>0
\end{equation} 
is simply the logarithm of the Ramanujan-Dedekind eta function $f(-q)$, where (see [7] chapters 21,22 and [3]):
$$
f(-q)=\prod^{\infty}_{n=1}\left(1-q^n\right)=\exp(-A)=
$$
\begin{equation}
=2^{1/3}\pi^{-1/2}q^{-1/24}(k_r)^{1/12}(k'_r)^{1/3}K(k_r)^{1/2}
\end{equation}
and (the definition of $E(x)$ will be used later): 
$$
K(x)=\frac{\pi}{2}\cdot {}_2F_1\left(\frac{1}{2},\frac{1}{2};1;x^2\right)\textrm{,  }E(x)=\frac{\pi}{2}\cdot {}_2F_1\left(-\frac{1}{2},\frac{1}{2};1;x^2\right),\eqno{(10.1)}
$$
with $0<k_r<1$ being the real solution of
$$
\frac{K\left(\sqrt{1-x^2}\right)}{K(x)}=\sqrt{r}.\eqno{(10.2)}
$$
Also 
\begin{equation}
\frac{dA}{dq}=\frac{e^{2x}}{4}\sum^{\infty}_{n=1}\frac{1}{\sinh^2(nx)}
\end{equation}
and
\begin{equation}
\sum^{\infty}_{n=1}\frac{1}{\sinh^2(nx)}=-4q\frac{d}{dq}\left(\log f(-q)\right)\textrm{, where }q=e^{-2x}\textrm{, }x>0,
\end{equation}
with
$$
\frac{dk_r}{dq}=\frac{-2k_r(k'_r)^2K(k_r)^2}{q\pi^2}\textrm{, }q=e^{-\pi\sqrt{r}}\textrm{, }r>0.\eqno{(11.1)}
$$   
Relation (11.1) is a result of Ramanujan (see [2] and [3]) for the first derivative of $k=k_r$ with respect to $q=e^{-\pi\sqrt{r}}$, $r>0$.\\ 
Also, (see [11]):
\begin{equation}
\sum^{\infty}_{n=1}\frac{(-1)^nn}{e^{2\pi n/\sqrt{r}}-1}=\frac{1}{8}-\frac{\sqrt{r}}{4\pi}+\frac{rK}{2\pi^2}(E-K).
\end{equation}
The functions $K=K(k_r)$ and $E=E(k_r)$ are called complete elliptic integrals of the first and second kind at singular values $k=k_r$. The singular modulus's $k_r$ are given from $K'/K=\sqrt{r}$ (this is (10.2), with $K'=K'(k_r)=K\left(k'_r\right)$ and $k'_r=\sqrt{1-k_r^2}$ is the complementary modulus) and help us give explicit evaluations, since: When $r$ is positive rational number, the values of $k_r$ are algebraic numbers. In this case also, the elliptic integrals $K=K(k_r),E=E(k_r)$ can evaluated and expressed in terms (as a finite product) of algebraic numbers, powers of $\pi$ and values of the Euler gamma $\Gamma(x)$ function (see [9],[10]).\\However, in what follows, we use the notation of $m(q)$ using relations (88) and (89) below. The notation $m(q)$ is more complete and more consisted to that of the elliptic singular modulus $k_r$, since $k_r=m\left(e^{-\pi\sqrt{r}}\right)$ and $k^{*}_r=m\left(-e^{-\pi\sqrt{r}}\right)$, $r>0$. Also $m(q)$ is defined in $q=e(z):=e^{2\pi i z}$, $Im(z)>0$.\\In the first three sections we use mainly the notation (10.2) but, one must be careful and have notation (88),(89) always in mind. In Section 4 and the rest sections, we assume mainly notation (88),(89) with $q=e(z)$, $z=x+i y$, $|x|<\frac{1}{2}$, $y>0$, to proceed further.

Continouing in view of [10] we can use elliptic alpha function to get a more detailed version of (13). It holds (for $q=e^{-\pi\sqrt{r}}$, $r>0$):
\begin{equation}
\alpha(r)=\frac{\pi}{4K^2}-\sqrt{r}\left(\frac{E}{K}-1\right).
\end{equation}
Solving with respect to $E$ the above formula we get
\begin{equation}
E=\frac{\pi}{4\sqrt{r}K}+K\left(1-\frac{\alpha(r)}{\sqrt{r}}\right).
\end{equation}
Hence (13) under the notes above becomes
\begin{equation}
2\sum^{\infty}_{n=1}\frac{n}{e^{4\pi n/\sqrt{r}}-1}-\sum_{n\geq 0,n-odd}\frac{n}{e^{2\pi n/\sqrt{r}}-1}=\frac{1}{8}-\frac{\sqrt{r}}{8\pi}-\alpha(r)\frac{\sqrt{r}K^2}{2\pi^2}.
\end{equation}   
Setting $r=4$ in (16) and using (7) along with the fact that $\alpha(4)=6-4\sqrt{2}$, we give the following example:
\begin{equation}
\sum_{n\geq0,n-odd}\frac{n}{e^{\pi n}-1}=-\frac{1}{24}+\frac{16\pi}{\Gamma\left(-\frac{1}{4}\right)^4}.
\end{equation}
Note also that, when $r$ is positive rational, then $\alpha(r)$ is algebraic number (see [10]).\\
Moreover in [14] we have proved\\
\\
\textbf{Theorem 2.}\\
Let $r>0$ and $q=e^{-\pi\sqrt{r}}$, then
\begin{equation}
1-24\sum^{\infty}_{n=1}\frac{n}{e^{\pi n\sqrt{r}}-1}
=\frac{6}{\pi\sqrt{r}}+\left(1+k^2_r-\frac{6 \alpha(r)}{\sqrt{r}}\right)\frac{4K^2}{\pi^2}.
\end{equation}
\\

A simple calculation using Theorem 2 and (see [10])
\begin{equation}
\alpha(4r)=(1+k_{4r})^2\alpha(r)-2\sqrt{r}k_{4r}
\end{equation}
can show us that:\\ 
\\
\textbf{Theorem 3.}\\
If $r>0$, then
\begin{equation}
1+24\sum_{n\geq 0,n-odd}\frac{n}{e^{\pi n\sqrt{r}}-1}=4\left(1+k_r^2\right)\frac{K^2}{\pi^2}.
\end{equation} 
\\
\textbf{Note.}\\
It is very interesting consequence of Theorems 2,3 and relation (16), the following modular property of elliptic alpha function 
\begin{equation}
\alpha\left(\frac{1}{r}\right)=\frac{1}{\sqrt{r}}-\frac{\alpha(r)}{r}\textrm{, }\forall r>0.
\end{equation} 
Relation (21) can be found in [10] chapter 5 pg.153.\\
\\

We continue with an evaluation theorem of Ramanujan (see [3] last chapter):\\
\\
\textbf{Theorem 4.}(Ramanujan)
\begin{equation}
1-24\sum^{\infty}_{n\geq0,n-odd}\frac{n}{e^{ny}+1}=z^2(1-2x),
\end{equation}
where $y=\pi\sqrt{r}=\pi K'/K$, $z=2K/\pi$ and $x=k^2=k_r^2$ (Ramanujan's notation).\\

Jacobi have given\\
\\
\textbf{Theorem 5.}(Jacobi)\\
For $r>0$ holds
\begin{equation}
\sum_{n\geq0,n-odd}\frac{1}{\cosh\left(n\pi\sqrt{r}/2\right)}=\frac{Kk_r}{\pi}.
\end{equation}
\\

Also, from relation (see [13])
\begin{equation}
\frac{\textrm{sn}(q,u)}{\textrm{cn}(q,u)\textrm{dn}(q,u)}=\frac{\pi}{2(k'_r)^2K}\tan\left(\frac{\pi u}{2K}\right)+\frac{2\pi}{(k'_r)^2K}\sum^{\infty}_{n=1}\frac{(-1)^nq^n}{1+q^n}\sin\left(\frac{n\pi u}{K}\right),
\end{equation}
we get\\
\\
\textbf{Theorem 6.}\\
If $r>0$, then
\begin{equation}
\sum^{\infty}_{n=0}\frac{(-1)^n}{e^{(2n+1)\pi\sqrt{r}}+1}=\frac{1}{4}-\frac{Kk'_r}{2\pi}.
\end{equation}

\section{Some General Properties and Related Series}

If $x>0$, then under some weak conditions on the sequence $X(n)$ we have the following formula for arithmetical functions $X(n)$:
\begin{equation}
\sum^{\infty}_{n=1}\frac{X(n)}{e^{nx}-1}=\sum^{\infty}_{n=1}e^{-nx}\sum_{d|n}X(d).
\end{equation}  
Hence the Ghost sum is the generating function of the divisor $d(n)=\sum_{d|n}1$ function i.e.\\
\\
\textbf{Proposition 1.}\\
\begin{equation}
g(x):=\sum^{\infty}_{n=1}\frac{1}{e^{nx}-1}=\sum^{\infty}_{n=1}d(n)q^n\textrm{, where }q=e^{-x}\textrm{, }x>0. 
\end{equation}
\\

This property makes it quite special. Another interesting thing with series of hyperbolic functions is that (as mentioned in the beginning of Introduction): ''almost all'' appear in the theory of Jacobian elliptic functions (see [8],[7]). Here is another example (we use the Ramanujan's notation) in [3] pg.176:\\
\\
\textbf{Proposition 2.}(Ramanujan),([3] pg.174)
\begin{equation}
\sec(\theta)+4\sum^{\infty}_{n=0}\frac{(-1)^n\cos((2n+1)\theta)}{e^{(2n+1)y}-1}=z\sec(\phi)\sqrt{1-x\sin^2(\phi)}.
\end{equation}

Relation (28) is quite close to what we search. 
Setting $\theta=0$ we get (in our notation):\\
\\
\textbf{Corollary 1.}\\
If $q=e^{-\pi\sqrt{r}}$, $r>0$
\begin{equation}
1+4\sum^{\infty}_{n=0}\frac{(-1)^n}{e^{(2n+1)\pi\sqrt{r}}-1}=\frac{2K}{\pi}.
\end{equation}
Hence also
\begin{equation}
1+4\sum_{n\geq0,n\equiv1(4)}\frac{1}{e^{n\pi\sqrt{r}}-1}-4\sum_{n\geq0,n\equiv3(4)}\frac{1}{e^{n\pi\sqrt{r}}-1}=\frac{2K}{\pi}.
\end{equation}
\\
Another result easy to check is\\
\\
\textbf{Theorem 7.}
\begin{equation}
\sum^{\infty}_{n=1}\frac{n}{\sinh^2(nx)}=-2\frac{d}{dx}\sum^{\infty}_{n=1}\frac{1}{e^{2nx}-1}\textrm{, }x>0.
\end{equation}
\\

The Fourier series of the Jacobi elliptic functions $\textrm{sn}$, $\textrm{cd}$ and $\textrm{cn}$, $\textrm{cd}$, $\textrm{sd}$ are:
$$
\textrm{sn}=\textrm{sn}(q,u)=\frac{2\pi}{Kk_r}\sum^{\infty}_{n=0}\frac{q^{n+1/2}\sin\left((2n+1)\frac{\pi}{2K}u\right)}{1-q^{2n+1}}
$$
$$
\textrm{cn}=\textrm{cn}(q,u)=\frac{2\pi}{Kk_r}\sum^{\infty}_{n=0}\frac{q^{n+1/2}\cos\left((2n+1)\frac{\pi}{2K}u\right)}{1+q^{2n+1}}
$$
$$
\textrm{cn}_1=\textrm{cn}_1(q,u)=\frac{2\pi}{Kk_r}\sum^{\infty}_{n=0}\frac{q^{n+1/2}\cos\left((2n+3)\frac{\pi}{2K}u\right)}{1+q^{2n+1}}
$$
$$
\textrm{sd}=\textrm{sd}(q,u)=\frac{2\pi}{Kk_rk'_r}\sum^{\infty}_{n=0}\frac{(-1)^nq^{n+1/2}\sin\left((2n+1)\frac{\pi}{2K}u\right)}{1+q^{2n+1}}
$$
$$
\textrm{cc}=\textrm{cc}(q,u)=\frac{2\pi}{Kk_r}\sum^{\infty}_{n=0}\frac{q^{n+1/2}\cos\left((2n+1)\frac{\pi}{2K}u\right)}{1+q^{2n-1}}=
$$
$$
=\frac{2\pi}{Kk_r}\frac{q^{1/2}\cos(z)}{1+q^{-1}}+q\cdot\textrm{cn}_1(q,u).
$$
$$
\textrm{cd}=\textrm{cd}(q,u)=\frac{2\pi}{Kk_r}\sum^{\infty}_{n=0}\frac{(-1)^nq^{n+1/2}\cos\left((2n+1)\frac{\pi}{2K}u\right)}{1-q^{2n+1}}
$$
$$
\textrm{dd}=\textrm{dd}(q,u)=\frac{2\pi}{Kk_r}\sum^{\infty}_{n=0}\frac{(-1)^nq^{n+1/2}\cos\left((2n+1)\frac{\pi}{2K}u\right)}{1-q^{2n-1}}
$$
and
\begin{equation}
\textrm{cd}_1=\textrm{cd}_1(q,u)=\frac{2\pi}{Kk_r}\sum^{\infty}_{n=0}\frac{(-1)^nq^{n+1/2}\cos\left((2n+3)\frac{\pi}{2K}u\right)}{1-q^{2n+1}}.
\end{equation}
The functions $\textrm{cn}_1$, $\textrm{cc}$, $\textrm{dd}$, $\textrm{cd}_1$ are not corresponding to elliptic functions and the notation is not the usual. For example $\textrm{cd}=\frac{cn}{dn}$, but $\textrm{dd}$ is not $\frac{dn}{dn}$ so one must be careful.\\

Beginning we can write
$$
\textrm{sn}(-q,u)=\frac{2\pi}{K^{*}k^{*}_r}\sum^{\infty}_{n=0}\frac{i(-1)^nq^{n+1/2}\sin\left((2n+1)\frac{\pi}{2K^{*}}u\right)}{1+q^{2n+1}}=
$$
$$
=i\frac{Kk_rk'_r}{K^{*}k^{*}_r}\frac{2\pi}{Kk_rk'_r}\sum^{\infty}_{n=0}\frac{(-1)^nq^{n+1/2}\sin\left((2n+1)\frac{\pi}{2K}u\frac{K}{K^{*}}\right)}{1+q^{2n+1}}=
$$
$$
=i\frac{Kk_rk'_r}{K^{*}k^{*}_r}\cdot\textrm{sd}\left(q,u\frac{K}{K^{*}}\right),\eqno{(32.1)}
$$
where we have set $K^{*}=K(k_{r_1})$ and $k^{*}=k_{r_1}$, with $r_1$ such that $e^{-\pi\sqrt{r_1}}=-e^{-\pi\sqrt{r}}$.\\ 
One can see that $k^2_rk^{{*}{2}}_r=k^2_r+k^{{*}{2}}_r$, or the equivalent $k'_r(k^{*}_r)'=1$. Hence from the modular identity
\begin{equation}
\frac{1}{\sqrt{1-x}}K\left(\frac{x}{x-1}\right)=K(x),
\end{equation}
we get
\begin{equation}
\frac{K(k_{r_1})}{K(k_{r})}=\frac{K^{*}}{K}=k'_r.
\end{equation}
\\
Hence a first result, that follows from (32.1), is the next:\\
\\
\textbf{Proposition 3.}\\
If $q=e^{-\pi\sqrt{r}}$, $r>0$, then
\begin{equation}
\textrm{sn}(-q,u)=k'_r\cdot\textrm{sd}\left(q,\frac{u}{k'_r}\right).
\end{equation}

Consider now Theorem 5 in the form
\begin{equation}
2\sum_{n\geq0,n-odd}\frac{q^{n/2}}{1+q^n}=\frac{Kk_r}{\pi}
\end{equation}
and set $q\rightarrow -q$. Then from the above relations (34) and $k^{*}_r=ik_r/k'_r$, we get\\
\\
\textbf{Theorem 8.}\\
If $q=e^{-\pi\sqrt{r}}$, $r>0$, then
\begin{equation}
2\sum^{\infty}_{n\geq0,n-odd}\frac{(-1)^{\frac{n-1}{2}}q^{n/2}}{1-q^n}=\sum^{\infty}_{n=0}\frac{(-1)^{n}}{\sinh\left((n+1/2)\pi\sqrt{r}\right)}=\frac{Kk_r}{\pi}.
\end{equation}
\\

Continuing if 
\begin{equation}
\textrm{ss}:=\textrm{ss}(q,u):=\frac{2\pi}{Kk_r}\sum^{\infty}_{n=0}\frac{q^{n+1/2}\sin\left((2n+1)\frac{\pi}{2K}u\right)}{1+q^{2n+1}},
\end{equation}
then from the above formula of $\textrm{cn}$ and the elementary trigonometric formula\\ $\cos(a+b)=\cos(a)\cos(b)-\sin(a)\sin(b)$, we get  
$$
\textrm{cc}(q,u)=\frac{2\pi}{Kk_r}\frac{q^{1/2}\cos(z)}{1+q^{-1}}+q\cdot\textrm{cn}_1(q,u)=\frac{2\pi}{Kk_r}\frac{q^{1/2}\cos(z)}{1+q^{-1}}+
$$
$$
+\frac{2\pi}{Kk_r}\sum^{\infty}_{n=0}\frac{q^{n+1/2+1}}{1+q^{2n+1}}\left(\cos((2n+1)z)\cos(2z)-\sin((2n+1)z)\sin(2z)\right)=
$$
$$
=\frac{2\pi}{Kk_r}\frac{q^{1/2}\cos(z)}{1+q^{-1}}+q\cos(2z) \textrm{cn}(q,u)-q\sin(2z)\textrm{ss}(q,u).
$$
Hence we get the next\\
\\
\textbf{Proposition 4.}\\
If $q=e^{-\pi\sqrt{r}}$, $r>0$
$$
\textrm{ss}(q,u)=\cot(2z)\textrm{cn}(q,u)-q^{-1}\csc(2z)\textrm{cc}(q,u)+\frac{\pi\sqrt{q}\csc(z)}{(1+q)k_rK}=
$$
\begin{equation}
=\textrm{cn}(q,u)\cot(2z)-\textrm{cn}_1(q,u)\csc(2z),
\end{equation}
where $z=\frac{\pi u}{2K}$.\\
\\ 
\textbf{Proposition 5.}\\
When $q=e^{-\pi\sqrt{r}}$, $r>0$, it holds
\begin{equation}
\textrm{cn}(-q,u)=\textrm{cd}\left(q,\frac{u}{k'_r}\right)
\end{equation}
\\
\textbf{Proof.}\\
From the $\textrm{cn}$ formula we have
$$
\textrm{cn}(-q,u)=\frac{2\pi}{K^{*}k^{*}_r}\sum^{\infty}_{n=0}\frac{(-1)^n q^{n+1/2}}{1-q^{2n+1}}\cos\left((2n+1)\frac{\pi u}{2K^{*}}\right)=
$$
$$
=i\frac{Kk_r}{K^{*}k^{*}_r}\frac{2\pi}{Kk_r}\sum^{\infty}_{n=0}\frac{(-1)^n q^{n+1/2}}{1-q^{2n+1}}\cos\left((2n+1)\frac{\pi u}{2K}\frac{K}{K^{*}}\right)=
$$
$$
=i\frac{Kk_r}{K^{*}k^{*}_r}\textrm{cd}\left(q,u\frac{K}{K^{*}}\right)=\textrm{cd}\left(q,\frac{u}{k'_r}\right).
$$
Since 
\begin{equation}
i\frac{Kk_r}{K^{*}k^{*}_r}=1.
\end{equation}
$qed$.

\section{The function $\textrm{cd}_1$: Evaluations and Properties}

\textbf{Theorem 9.}\\
If $A=iq^{1/2}e^{i\pi u/(2K)}$ and $u\in \textbf{C}$, with $|q^{1/2}e^{i\pi u/(2K)}|<1$, then
\begin{equation}
\frac{2\pi }{Kk_r}\sum^{\infty}_{n=0}\frac{A^{2n+1}}{1-q^{2n+1}}=-\textrm{cd}\left(q,u\right)\cot\left(\frac{u\pi}{K}\right)+\textrm{cd}_1\left(q,u\right)\csc\left(\frac{u\pi}{K}\right)+i\cdot\textrm{cd}\left(q,u\right).
\end{equation}
\\
\textbf{Proof.}\\
It is easy to see someone that
$$
\frac{2\pi }{Kk_r}\sum^{\infty}_{n=0}\frac{A^{2n+1}}{1-q^{2n+1}}=
i\cdot\textrm{cn}(-q,u k'_r)-\textrm{ss}\left(-q,u k'_r\right)=
$$
$$
=i\cdot \textrm{cd}\left(q,u\right)-\textrm{ss}\left(-q,u k'_r\right).\eqno{(42.1)}
$$
But from Propositions 4,5 we have
$$
\textrm{ss}(-q,u)=\textrm{cd}\left(q,\frac{u}{k'_r}\right)\cot(2z^{*})-\textrm{cn}_1(-q,u)\csc(2z^{*})\textrm{, }z^{*}=\frac{\pi u}{2K^{*}}.
$$
Also
\begin{equation}
\textrm{cn}_1(-q,u)=\textrm{cd}_1\left(q,\frac{u}{k'_r}\right).
\end{equation}
Hence
\begin{equation}
\textrm{ss}(-q,u)=\textrm{cd}\left(q,\frac{u}{k'_r}\right)\cot\left(\frac{2z}{k'_r}\right)-\textrm{cd}_1\left(q,\frac{u}{k'_r}\right)\csc\left(\frac{2z}{k'_r}\right).
\end{equation}
From the above relations we get the proof. $qed$\\
\\
\textbf{Theorem 10.}\\
If $0<|\lambda|<1$, $r>0$, then
$$
\frac{\pi}{Kk_r}\sum^{\infty}_{n=0}\frac{(-1)^ne^{-\pi\sqrt{r}(n+1/2)\lambda}}{\sinh\left((n+1/2)\pi\sqrt{r}\right)}=\textrm{cd}\left(q,\lambda i K'\right)\coth\left(\lambda\pi\sqrt{r}\right)-
$$
\begin{equation}
-\textrm{cd}_1\left(q,\lambda i K'\right)\textrm{csch}\left(\lambda\pi\sqrt{r}\right)+\textrm{cd}\left(q,\lambda i K'\right).
\end{equation}
\\
\textbf{Proof.}\\
Set $u=\lambda i K'$ in relation (42) of Theorem 9. $qed$\\
\\
\textbf{Theorem 11.}\\
If $q=e^{-\pi\sqrt{r}}$, $r>0$ and $\nu\in \textbf{C}^{*}$, such that $2/\nu$ not integer, then
$$
\frac{2\pi}{Kk_r}\sum^{\infty}_{n=0}\frac{q^{(2n+1)(1/2+1/\nu)}}{1-q^{2n+1}}=i\cdot\textrm{sn}\left(q,\frac{2iK'}{\nu}\right)\coth\left(\frac{2\pi\sqrt{r}}{\nu}\right)+
$$
\begin{equation}
+i\cdot\textrm{cd}_1\left(q,-K+\frac{2iK'}{\nu}\right)\textrm{csch}\left(\frac{2\pi\sqrt{r}}{\nu}\right)+i\cdot\textrm{sn}\left(q,\frac{2iK'}{\nu}\right).
\end{equation}
In case that $r\in\textbf{Q}^{*}_{+}$ and $\nu\in\textbf{Q}^{*}_{+}-\textbf{Z}$, then $\textrm{sn}\left(q,\frac{2iK'}{\nu}\right)$ is algebraic number.\\ 
\\
\textbf{Proof.}\\
If we replace 
$$
u=u_1:=-i\log\left(-iq^{1/\nu}\right)\frac{2K}{\pi}=\left(-\frac{\pi}{2}+i\frac{\pi\sqrt{r}}{\nu}\right)\frac{2K}{\pi}=
$$
\begin{equation}
=\left(-1+i\frac{2\sqrt{r}}{\nu}\right)K=-K+\frac{2iK'}{\nu},
\end{equation}
then we will have $A=q^{1/2+1/\nu}$. From the relations (see [8]):
\begin{equation}
\textrm{cn}\left(q,w+K\right)=-k'_r\textrm{sn}(q,w)/\textrm{dn}(q,w)
\end{equation}
and
\begin{equation}
\textrm{dn}\left(q,w+K\right)=k'_r/\textrm{dn}(q,w)
\end{equation}
we have
\begin{equation}
\textrm{cd}\left(q,w+K\right)=-\textrm{sn}(q,w)
\end{equation}
and easy
$$
\textrm{cd}\left(q,w-K\right)=-\textrm{sn}\left(q,w-2K\right)=\textrm{sn}(q,w),
$$
since $\textrm{sn}\left(q,w+2K\right)=-\textrm{sn}(q,w)$.\\Also easily we get
$$
\csc\left(\frac{u_1\pi}{K}\right)=i\cdot\textrm{csch}\left(\frac{2\pi\sqrt{r}}{\nu}\right)
$$ 
and 
$$
\cot\left(\frac{u_1\pi}{K}\right)=-i\cdot\textrm{coth}\left(\frac{2\pi\sqrt{r}}{\nu}\right).
$$
From the above and Theorem 9 we get the result. $qed$\\ 
\\
\textbf{Corollary 2.}\\
If $q=e^{-\pi\sqrt{r}}$, $r>0$, then\\
1)
\begin{equation}
\lim_{y\rightarrow K}\frac{\textrm{cd}_1\left(q,y\right)}{y-K}=1+\frac{2\pi^2}{K^2k_r}\sum_{n\geq0\textrm{, }n-odd}\frac{q^{n/2}}{1-q^n},
\end{equation}
2) $\textrm{cd}_1(q,0)=1$, $\textrm{cd}_1(q,K)=0$, $\textrm{cd}_1(q,2K)=-1$ and 
$$
\textrm{cd}_1(q,u+2K)=-\textrm{cd}_1(q,u)\textrm{, }\forall u\in\textbf{R}.
$$
\\
\textbf{Proof.}\\
Taking the limit $\nu\rightarrow\infty$ in (46) and using
$$
\lim_{\nu\rightarrow\infty}\textrm{sn}\left(q,\frac{2iK'}{\nu}\right)=0,
$$
$$
\lim_{\nu\rightarrow\infty}\textrm{sn}\left(q,\frac{2iK'}{\nu}\right)\coth\left(\frac{2\pi\sqrt{r}}{\nu}\right)=\frac{iK'}{\pi\sqrt{r}},
$$
we get easily the result. $qed$\\
\\
\textbf{Theorem 12.}\\
If $q=e^{-\pi\sqrt{r}}$, $r>0$, then
\begin{equation}
\textrm{cd}_1(q,iK')=\frac{1}{qk_r}-\frac{\sinh(\pi\sqrt{r})}{k_r}\left(1-\frac{\pi}{2K}\right).
\end{equation}
\\
\textbf{Proof.}\\
Set $u=iK'$ in (42). Then using relation $\textrm{cd}\left(q,iK'\right)=\frac{1}{k_r}$ (see [8]), we get 
$$
\textrm{cd}_1\left(q,iK'\right)=\frac{1}{qk_r}-\frac{2\pi\sinh\left(\pi\sqrt{r}\right)}{Kk_r}\sum^{\infty}_{n=0}\frac{(-1)^n}{e^{(2n+1)\pi\sqrt{r}}-1}.
$$
The result follows from Corollary 1. $qed$\\
\\
\textbf{Notes.}\\
i) Numerical values of $\textrm{cd}_1(q,iK')$ can given using Theorems 13 and 14 bellow.\\ 
ii) Formula (32) does not converge for these values. So we can speak here for analytic continuation of $\textrm{cd}_1$.\\
\\
\textbf{Corollary 3.}\\
If $q=e^{-\pi\sqrt{r}}$, $r>0$, then
\begin{equation}
\textrm{cd}_1\left(q,\frac{iK'}{2}\right)=\frac{1}{\sqrt{qk_r}}-\frac{\pi\sinh\left(\frac{\pi\sqrt{r}}{2}\right)}{Kk_r}\sum^{\infty}_{n=0}\frac{(-1)^ne^{-(n+1/2)\pi\sqrt{r}/2}}{\sinh\left((n+1/2)\pi\sqrt{r}\right)}
\end{equation}
\\
\textbf{Corollary 4.}\\
If $q=e^{-\pi\sqrt{r}}$, $r>0$, then
\begin{equation}
\textrm{cd}_1\left(q,\frac{K}{2}\right)=\frac{1}{\sqrt{1+k'_r}}-\frac{\pi\sqrt{8}}{Kk_r}\sum_{n\geq0\textrm{,  }n\equiv 1(4)}\frac{(-1)^{\frac{n-1}{4}}q^{n/2}}{1-q^n}
\end{equation}
and
\begin{equation}
-\sum^{\infty}_{n=1}\chi(n)\frac{q^{n/2}}{1-q^n}=\frac{Kk_r}{\sqrt{2}\pi\sqrt{1+k'_r}},
\end{equation}
where $\chi(n)=\left(\frac{n+2}{8}\right)$ and $\left(\frac{n}{m}\right)$ is the usual Jacobi symbol.\\ 
\\
\textbf{Proof.}\\
Setting $u=\frac{K}{2}$ in Theorem 9, we get $a=iq^{1/2}e^{\pi i/4}=q^{1/2}e^{3\pi i/4}$, $\cot\left(\frac{\pi}{2}\right)=0$, $\csc\left(\frac{\pi}{2}\right)=1$. Also it holds (see [8]):
\begin{equation}
\textrm{cn}\left(q,\frac{K}{2}\right)=\frac{\sqrt{k'_r}}{\sqrt{1+k'_r}}\textrm{ and  }\textrm{dn}\left(q,\frac{K}{2}\right)=\sqrt{k'_r}.
\end{equation}
Hence
\begin{equation}
\textrm{cd}\left(q,\frac{K}{2}\right)=\frac{\textrm{cn}\left(q,\frac{K}{2}\right)}{\textrm{dn}\left(q,\frac{K}{2}\right)}=\frac{1}{\sqrt{1+k'_r}}.
\end{equation}
From the above we can evaluate
\begin{equation}
\textrm{cd}_1\left(q,\frac{K}{2}\right)=e^{3\pi i/4}\frac{2\pi}{Kk_r}\sum^{\infty}_{n=0}\frac{q^{n+1/2}e^{3\pi i n/2}}{1-q^{2n+1}}-i\cdot\textrm{cd}\left(q,\frac{K}{2}\right).
\end{equation}
Hence taking the real and imaginary parts of the above equation we deduce 
\begin{equation}
\sum_{n\geq 0\textrm{, }n\equiv 1(4)}\frac{(-1)^{\frac{n-1}{4}}q^{n/2}}{1-q^n}+\sum_{n\geq 0\textrm{, }n\equiv 3(4)}\frac{(-1)^{\frac{n-3}{4}}q^{n/2}}{1-q^n}=\frac{\sqrt{2}}{\sqrt{1+k'_r}}\frac{Kk_r}{2\pi}
\end{equation}
and
\begin{equation}
\frac{\pi\sqrt{2}}{Kk_r}\sum_{n\geq 0\textrm{, }n\equiv 3(4)}\frac{(-1)^{\frac{n-3}{4}}q^{n/2}}{1-q^n}-\frac{\pi\sqrt{2}}{Kk_r}\sum_{n\geq 0\textrm{, }n\equiv 1(4)}\frac{(-1)^{\frac{n-1}{4}}q^{n/2}}{1-q^n}=\textrm{cd}_1\left(q,\frac{K}{2}\right).
\end{equation}
$qed$\\

\section{The $\theta(q,u)$ modular angle and the functions $\textrm{ss}$ and $\textrm{cd}_1$}

In this section we give rise to the evaluation of $\int \textrm{ss}(q,u)du$ and hence to $\textrm{ss}(q,u)$ and  $\textrm{cd}_1(q,u)$, up to a Jacobi modular angle. After that we proceed further with the analytical study of the functions $\textrm{ss}$, $\textrm{cd}_1$. We also give some applications concerning the area of  Jacobian elliptic functions, hyperbolic series and continued fractions. These connections are indeed very interesting since they fill large gaps in the literature. We start with a continued fraction expansion of Ramanujan.\\

Ramanujan has stated that (see [3] pg.21):\\
\\
\textbf{Theorem 13.}\\
Suppose that $q$, $a$ and $b$ are complex numbers with $\left|q\right|<1$, or that $q,a$, and $b$ are complex numbers with $a=bq^m$ for some integer $m$. Then
$$
U=U(a,b;q)=\frac{(-a;q)_{\infty}(b;q)_{\infty}-(a;q)_{\infty}(-b;q)_{\infty}}{(-a;q)_{\infty}(b;q)_{\infty}+(a;q)_{\infty}(-b;q)_{\infty}}=
$$  
\begin{equation}
=\frac{a-b}{1-q+}\frac{(a-bq)(aq-b)}{1-q^3+}\frac{q(a-bq^2)(aq^2-b)}{1-q^5+}\frac{q^2(a-bq^3)(aq^3-b)}{1-q^7+}\ldots.
\end{equation}
\\

Also holds (see [12])
\begin{equation}
\left(-1+\frac{2}{1-U(a,b;q)}\right)^2=\frac{\left(-1+\frac{2}{1-u_0(q,a)}\right)}{\left(-1+\frac{2}{1-u_0(q,b)}\right)},
\end{equation}
where the $u_0(q,a)$ are introduced as below:\\

We define 
\begin{equation}
(z;q)_{\infty}:=\prod^{\infty}_{n=0}(1-zq^n)\textrm{, }|q|<1\textrm{, }z\in\textbf{C}.
\end{equation}
Also
\begin{equation}
P:=\left(\frac{(-A;q)_{\infty}}{(A;q)_{\infty}}\right)^2.
\end{equation}
Then
\begin{equation}
u_0(q,A):=\frac{P-1}{P+1}.
\end{equation}
Then we have (see [12]):
\begin{equation}
\log\left(-1+\frac{2}{1-u_0(q,A)}\right)=\log P
\end{equation}
and
\begin{equation} \log P=\log\left(-1+\frac{2}{1-u_0(q,A)}\right)=4\sum^{\infty}_{n=0}\frac{A^{2n+1}}{(2n+1)(1-q^{2n+1})}.
\end{equation}
Let $u$ be a complex number. We set $A=iq^{1/2}e^{i\pi t/(2K)}$ in (67) and derivate with respect to parameter $t$. Then setting the value $t=u$, we get
$$
\left[\frac{d}{dt}\log\left(-1+\frac{2}{1-u_0(q,A)}\right)\right]_{t=u}=4\left[\frac{d}{dt}\sum^{\infty}_{n=0}\frac{A^{2n+1}}{(2n+1)(1-q^{2n+1})}\right]_{t=u}=
$$
\begin{equation}
=\frac{2\pi i}{K}\sum^{\infty}_{n=0}\frac{A^{2n+1}}{1-q^{2n+1}}.
\end{equation}
Hence using relation (42.1) we get
\begin{equation}
\frac{d}{dt}\log\left(-1+\frac{2}{1-u_0(q,A)}\right)
=-k_r\textrm{cd}(q,t)-i\cdot k_r\textrm{ss}\left(-q,  k'_r t\right),
\end{equation}
where
\begin{equation}
\textrm{ss}=\textrm{ss}(q,u):=\frac{2\pi}{Kk_r}\sum^{\infty}_{n=0}\frac{q^{n+1/2}\sin\left((2n+1)\frac{\pi}{2K}u\right)}{1+q^{2n+1}}.
\end{equation}
After integration of (69) with respect to $t$ from $0$ to $u$, we get the next:\\
\\
\textbf{Theorem 14.}\\
If $A=iq^{1/2}e^{i\pi u/(2K)}$, then
$$
\log\left(-1+\frac{2}{1-u_0(q,A)}\right)=-\log\left(\textrm{nd}(q,u)+k_r\textrm{sd}\left(q,u\right)\right)+
$$
\begin{equation}
+4i\sum^{\infty}_{n=0}\frac{(-1)^nq^{n+1/2}\cos\left((2n+1)\frac{\pi}{2K}u\right)}{(2n+1)\left(1-q^{2n+1}\right)}.
\end{equation}
\\
\textbf{Proof.}\\
We give more details for how we arived into (71). We have the relation
$$
\textrm{ss}(-q,u)=\frac{2\pi i}{K^{*}k_r^{*}}\sum^{\infty}_{n=0}\frac{(-1)^nq^{n+1/2}\sin\left((2n+1)\frac{\pi }{2K^{*}}u\right)}{1-q^{2n+1}},
$$
where $K^{*}=K(k^{*}_r)$ and $k^{*}_r$ corresponds to the change of sign $q\rightarrow -q$. Hence from (34) and the related identities   
$$
k^{*}_r=\frac{ik_r}{k'_r}\eqno{(71.0)}
$$
and
$$
K^{*}=k'_rK\textrm{, where }k'_r=\sqrt{1-k_r^2},
$$
we get
$$
\textrm{ss}(-q,uk'_r)=\frac{2\pi}{Kk_r}\sum^{\infty}_{n=0}\frac{(-1)^nq^{n+1/2}\sin\left((2n+1)\frac{\pi}{2K}u\right)}{1-q^{2n+1}}.\eqno{(71.1)}
$$
Hence integrating (69) we get (71) using (see [8])
\begin{equation}
\int \textrm{cd}(q,t)dt=\log\left(\textrm{nd}(q,t)+k_r\textrm{sd}(q,t)\right)
\end{equation}
and (71.1). $qed$\\
\\
\textbf{Theorem 15.}\\
If $q=e^{-\pi\sqrt{r}}$, $r>0$ and $A=iq^{1/2}e^{i\pi u/(2K)}$, $u$ real number, then 
\begin{equation}
Re\left[\log\left(-1+\frac{2}{1-u_0(q,A)}\right)\right]
=-\log\left(\textrm{nd}(q,u)+k_r\textrm{sd}(q,u)\right)
\end{equation}
and
\begin{equation}
\textrm{cd}_1(q,u)=\textrm{cd}(q,u)\cos\left(\frac{\pi u}{K}\right)+2k_r^{-1}\sin\left(\frac{\pi u}{K}\right)Im\left[\frac{d}{d u}\log\left(\frac{\left(-A;q\right)_{\infty}}{\left(A;q\right)_{\infty}}\right)\right].
\end{equation}
\\
\textbf{Proof.}\\
See the proof of Theorem 16 below.\\
\\

From [12] we have the next continued fraction expansion:
\begin{equation}
u_0(q,a)=\frac{2a}{1-q+}\frac{a^2(1+q)^2}{1-q^3+}\frac{a^2q(1+q^2)^2}{1-q^5+}\frac{a^2q^2(1+q^3)^2}{1-q^7+}\ldots\textrm{, }|q|,|q/a|<1.
\end{equation}
This continued fraction can be used to get numerical and arithmetical evaluations (verifications) for complex values of $\textrm{cd}_1(q,u)$, $u\in\textbf{C}$ (such $\textrm{cd}_1(q,i\cdot nK')\textrm{, } n\in\textbf{Q}^{*}_{+}$). More precisely we have the next analytic continuation theorem:\\ 
\\
\textbf{Theorem 16.}\\
If $q=e^{-\pi\sqrt{r}}$, $r>0$ and $A=iq^{1/2}e^{i\pi t/(2K)}$, $t$ parameter and $u\in\textbf{C}$ such that $|q^{1/2}e^{i\pi u/(2K)}|<1$, then
$$
\textrm{cd}_1(q,u)=\textrm{cd}(q,u)\cos\left(\frac{\pi u}{K}\right)-i\cdot\textrm{cd}(q,u)\sin\left(\frac{\pi u}{K}\right)-i\cdot k_r^{-1}\sin\left(\frac{\pi u}{K}\right)\times
$$
\begin{equation}
\times\left[\frac{d}{dt}\log\left(-1+\frac{2}{1+}\frac{-2A}{1-q+}\frac{A^2(1+q)^2}{1-q^3+}\frac{A^2q(1+q^2)^2}{1-q^5+}\frac{A^2q^2(1+q^3)^2}{1-q^7+}\ldots\right)\right]_{t=u}
\end{equation}
and
$$
\textrm{cd}_1(q,u)=\textrm{cd}(q,u)\cos\left(\frac{\pi u}{K}\right)-i\cdot\textrm{cd}(q,u)\sin\left(\frac{\pi u}{K}\right)
-
$$
\begin{equation}
-2i\cdot k_r^{-1}\sin\left(\frac{\pi u}{K}\right)\left[\frac{d}{dt}\log\left(\frac{(-A;q)_{\infty}}{(A;q)_{\infty}}\right)\right]_{t=u}.
\end{equation}
\\
\textbf{Proof.}\\
Using (64),(65),(66),(67),(68),(75) and (42), we get the two results. $qed$\\
\\
\textbf{Theorem 17.}\\
Let $1/\nu$ be positive integer. Let also $m$ be even integer and $\nu_1=2/\nu$, then if $q=e^{-\pi\sqrt{r}}$, $r>0$ we have
$$
\textrm{cd}_1\left(q,mK+\nu_1iK'\right)=(-1)^{m/2}e^{\nu_1\pi \sqrt{r}}-
$$
\begin{equation}
-(-1)^{m/2}\sinh\left(\nu_1\pi\sqrt{r}\right)\left(1-\frac{2\pi }{Kk_r}\sum^{1/\nu-1}_{j=0}\frac{q^{j+1/2}}{1+q^{2j+1}}
\right).
\end{equation} 
\\
\textbf{Proof.}\\
The proof follows from Theorems 9,11 and the identities (67),(68) along with
\begin{equation}
\sum^{\infty}_{n=0}\frac{(-1)^nq^{(2n+1)(l+1/2)}}{1-q^{2n+1}}=-\sum^{l-1}_{j=0}\frac{q^{j+1/2}}{1+q^{2j+1}}+\frac{Kk_r}{2\pi}\textrm{, }l\in\textbf{N}\textrm{, }q=e^{-\pi\sqrt{r}}\textrm{, }r>0,
\end{equation}
and 
\begin{equation}
\textrm{cd}(q,mK+i\nu_1K')=(-1)^{m/2},
\end{equation}
with $m$ even integer and $\nu_1=2/\nu$ even positive integer. $qed$\\
\\

If $\theta(a)$ and $\theta(b)$ denote modular angles, where $\theta(x)=\theta(q,x)$ (see [7]) such that
\begin{equation}
e^{-\theta(q,x)}=\frac{\prod^{\infty}_{n=0}\left(1-q^{n+x}\right)}{\prod^{\infty}_{n=0}\left(1+q^{n+x}\right)}\textrm{, }x>0,
\end{equation}
then
\begin{equation}
\theta(q,x)=2\sum^{\infty}_{n=0}\textrm{arctanh}\left(q^{n+x}\right).
\end{equation}
Then also holds (see [12]):
$$
2\theta(q,a)-2\theta(q,b)=
$$
\begin{equation}
=2\sum^{\infty}_{n=0}\textrm{arctanh}\left(q^{n+a}\right)-2\sum^{\infty}_{n=0}\textrm{arctanh}\left(q^{n+b}\right)=\log\left(-1+\frac{2}{1-U(q^a,q^b,q)}\right).
\end{equation}

A first result that one can derive immediately from the definition of the modular angles is of course:\\
\\
\textbf{Theorem 17.1}
$$
\log\left(-1+\frac{2}{1-u_0\left(q,q^a\right)}\right)=2\theta(q,a)\textrm{, }|q|<1\textrm{, }a>0\eqno{(83.1)}
$$
\\
\textbf{Notes.}\\
i) In case $a,b\in \textbf{R}$ and $\{a\}=\{b\}=c$, ($\{x\}$ is the fractional part of $x$) we have evaluated successfully (83) in [12].\\ 
ii) Also in [12] relation (73) have been proved and function $\textrm{ss}$ has been considered.\\ 
iii) In what it follows we always denote with $\theta$ the modular angle $\theta(x)=\theta(q,x)$. Similar notations such $\theta_0$,$\theta_1$,$\theta_2$ of Theorem 18 below have nothing to do with the above notation of modular angle.\\
iv) From the identities ($q=e^{-\pi\sqrt{r}}$, $r>0$):
\begin{equation}
f(-q):=\prod^{\infty}_{n=1}(1-q^n)=2^{1/3}\pi^{-1/2}q^{-1/24}k_r^{1/12}(k'_r)^{1/3}K^{1/2},
\end{equation}
\begin{equation}
\prod^{\infty}_{n=1}(1+q^n)=2^{-1/6}q^{-1/24}k_r^{1/12}(k'_r)^{-1/6},
\end{equation}
and relation (81), if we assume that $a\in\textbf{N}=\{1,2,3,\dots\}$, then we easily can obtain the next special values
\begin{equation}
e^{-\theta(q,a)}=\exp\left[-2\cdot\textrm{arctanh}\left(q^a\right)\right]\sqrt{\frac{2Kk'_r}{\pi}}\prod^{a}_{j=1}\frac{1+q^j}{1-q^j}\textrm{, }q=e^{-\pi\sqrt{r}}\textrm{, }r>0.
\end{equation}
v) We will also use the generalization $m(q)$ of $k_r$ as defined in (88),(89) below. For the restricted case $q=e^{-\pi\sqrt{r}}$, $r>0$, we have $m\left(q\right)=k_r$, $r>0$.\\
\\
\textbf{Theorem 18.}\\
If $a,x,y$ are real with $0<a<1$, $x+\frac{1}{2}\in\left(-\frac{1}{2},\frac{1}{2}\right)$ and $y>0$, $q=e(z):=\exp(2\pi i z)$, where $z=x+iy$, and $\theta_1=-K+(4a-2)xK$,  $\theta_2=(4a-2)yK$, $\theta_0=\theta_1+i\theta_2$, $K=K(m(q))$, then
$$
q^a=iq^{1/2}e^{i\pi \theta_0/(2K)}\eqno{(86.1)}
$$
and
$$
\sum^{\infty}_{n=0}\frac{(-1)^nq^{n+1/2}}{(2n+1)(1-q^{2n+1})}\cos\left((2n+1)\frac{\pi}{2K}\theta_0\right)=
$$
\begin{equation}
=\frac{1}{2i}\theta(q,a)+\frac{1}{4i}\log\left(\textrm{nd}(q,\theta_0)+m(q)\cdot\textrm{sd}(q,\theta_0)\right).
\end{equation}
\\
\textbf{Proof.}\\
If $a,x,y$ and $\theta_1,\theta_2$ are as in the statement of the theorem and $K=K(m(q))$, $K'=K(m'(q))$, with $q=\exp(2\pi i z)$, $z=x+iy$, $x+\frac{1}{2}\in\left(-\frac{1}{2},\frac{1}{2}\right)$, $y>0$, then we define  
\begin{equation}
m(q):=\frac{\vartheta_2(q)^2}{\vartheta_3(q)^2}\textrm{, }m'(q):=\sqrt{1-m(q)^2},
\end{equation} 
where
\begin{equation}
\vartheta_2(q):=\sum^{\infty}_{n=-\infty}q^{\left(n+1/2\right)^2}\textrm{, }\vartheta_3(q):=\sum^{\infty}_{n=-\infty}q^{n^2}\textrm{, }|q|<1.
\end{equation}
For this definition it holds
$$
m(-q)=i\frac{m(q)}{m'(q)}\eqno{(89.1)}
$$
and
$$
i\frac{K(m'(q))}{K(m(q))}=2z.\eqno{(89.2)}
$$
Also then $q^a=iq^{1/2}e^{\pi i \theta_0/(2K)}$ is satisfied. Hence from Theorems 14,17.1 we get the result. $qed$\\
\\

Integration of (71.1) lead us to 
$$
\int^{\theta_0}_{0} \textrm{ss}\left(-q,m'(q) t\right)dt=-\frac{4}{m(q)}\sum^{\infty}_{n=0}\frac{(-1)^nq^{n+1/2}\cos\left((2n+1)\frac{\pi}{2K}\theta_0\right)}{(2n+1)\left(1-q^{2n+1}\right)}+
$$
\begin{equation}
+\frac{4}{m(q)}\sum^{\infty}_{n=0}\frac{(-1)^nq^{n+1/2}}{(2n+1)\left(1-q^{2n+1}\right)}.
\end{equation}
Hence more clear:\\
\\
\textbf{Theorem 19.}\\
Let $a,x,y,\theta_0,\theta_1,\theta_2$ be as in Theorem 18, then the nome is $q^a=iq^{1/2}e^{\pi i \theta_0/(2K)}$ and hold the following relations:
$$
\int^{\theta_0}_{0}\textrm{ss}(-q,m'(q) t)dt=\frac{i}{m(q)}\log\left(e^{2\theta(q,a)}\left(\textrm{nd}(q,\theta_0)+m(q)\cdot\textrm{sd}(q,\theta_0)\right)\right)+
$$
\begin{equation}
+\frac{4}{m(q)}\sum^{\infty}_{n=0}\frac{(-1)^nq^{n+1/2}}{(2n+1)\left(1-q^{2n+1}\right)}
\end{equation}
and
\begin{equation}
\textrm{ss}\left(-q,m'(q)\theta_0\right)=-\frac{2\pi}{m(q)K}\sum^{\infty}_{n=0}\frac{q^{a(2n+1)}}{1-q^{2n+1}}+i\cdot\textrm{cd}(q,\theta_0).
\end{equation}
\\
\\
\textbf{Proof.}\\
Relation (91) follows from (90) and Theorem 18. The relation (92) follows from (69),(68) and the nome $q^a=iq^{1/2}e^{\pi i \theta_0/(2K)}$. $qed$\\
\\
\textbf{Theorem 20.}\\
If $q=e(z)$, $z=x+iy$ and  $x\in\left(-\frac{1}{2},\frac{1}{2}\right)$, $y>0$, then
$$
4\pi i z\sum^{\infty}_{n=0}\frac{(-1)^nq^{n+1/2}}{1-q^{2n+1}}\sin\left((2n+1)\frac{\pi\theta_0}{2K}\right)=
$$
\begin{equation}
=-\left(\frac{\partial\theta(q,u)}{\partial u}\right)_{u=a}-2zKm(q)\textrm{cd}\left(q,\theta_0\right).
\end{equation}
If $x+\frac{1}{2}\in\left(-\frac{1}{2},\frac{1}{2}\right)$, then
$$
4\pi i \left(z+\frac{1}{2}\right)\sum^{\infty}_{n=0}\frac{q^{n+1/2}}{1+q^{2n+1}}\sin\left((2n+1)\frac{\pi\theta'_0}{2K}\right)
=i\left(\frac{\partial\theta(-q,u)}{\partial u}\right)_{u=a}-
$$
\begin{equation}
-2\left(z+\frac{1}{2}\right)Km(q)\textrm{cn}\left(q,\theta'_0\right),
\end{equation}
where $\theta'_0=-K+\left(4a-2\right) \left(z+\frac{1}{2}\right)K$.\\ 
\\
\textbf{Proof.}\\
For to prove (93), differentiate (87) with respect to $a$ using $\frac{d\theta_0}{da}=4zK$ and
$$
\frac{d}{d\theta_0}\log\left(\textrm{nd}(q,\theta_0)+m(q)\cdot \textrm{sd}(q,\theta_0)\right)=m(q)\textrm{cd}(q,\theta_0).
$$ 
For to prove relation (94), set in (87) where $q\rightarrow -q$ and do the same.\\
\\

This also implies the following evaluation of $\textrm{ss}(q,u)$ function:\\
\\
\textbf{Theorem 20.1}\\
If $q=e(z)$, $z=x+iy$, $y>0$, $x+\frac{1}{2}\in\left(-\frac{1}{2},\frac{1}{2}\right)$ and $\theta'_0=-K+(4a-2)(z+\frac{1}{2})K$, then
\begin{equation}
\textrm{ss}\left(q,\theta'_0\right)
=\frac{1}{(2z+1)m(q)K}\left(\frac{\partial\theta^{*}(q,u)}{\partial u}\right)_{u=a}+i\cdot\textrm{cn}\left(q,\theta'_0\right),
\end{equation}
where $\theta^{*}(q,u):=\theta(-q,u)$.\\
\\

We introduce the notation
$$
\frac{d\theta(a)}{da}:=\left(\frac{\partial\theta(q,u)}{\partial u}\right)_{u=a}\textrm{ and }\frac{d\theta^{*}(a)}{da}:=\left(\frac{\partial\theta(-q,u)}{\partial u}\right)_{u=a}.\eqno{(95.1)}
$$

If the nome is $q^a=iq^{1/2}e^{i\pi\theta_0/(2K)}$, then from (68) and (83.1) we have
\begin{equation}
\frac{\pi i}{K}\sum^{\infty}_{n=0}\frac{q^{a(2n+1)}}{1-q^{2n+1}}=\frac{d\theta(a)}{d\theta_0}.
\end{equation}   
\\
Relation (96) can also written as
$$
\frac{d\theta_0}{da}\cdot\sum^{\infty}_{n=0}\frac{q^{a(2n+1)}}{1-q^{2n+1}}=-\frac{iK}{\pi}\frac{d\theta(a)}{da}.
$$
Hence we get the next\\
\\
\textbf{Theorem 20.2}\\
Let $q=e(z)$, $z=x+iy$, $y>0$, $x+\frac{1}{2}\in\left(-\frac{1}{2},\frac{1}{2}\right)$ and $\theta_0$ as in Theorem 18 i.e it holds $q^a=iq^{1/2}e^{i\pi\theta_0/(2K)}$, then
\begin{equation}
\frac{d\theta_0}{da}=4zK=2iK'\textrm{, }a>0
\end{equation} 
and
\begin{equation}
\sum^{\infty}_{n=0}\frac{q^{a(2n+1)}}{1-q^{2n+1}}=\frac{1}{4\pi i z}\frac{d\theta(a)}{da}=\frac{1}{4\pi i z}\left(\frac{\partial\theta(q,u)}{\partial u}\right)_{u=a}\textrm{, }a>0.
\end{equation}
\\

Continuing our arguments, from relation (44) and  definition of $\theta_0$ we have:
$$
\textrm{ss}(-q,m'(q)\theta_0)=\textrm{cd}(q,\theta_0)\cot\left(\frac{\pi\theta_0}{K}\right)-\textrm{cd}_1(q,\theta_0)\csc\left(\frac{\pi\theta_0}{K}\right),
$$
where $\textrm{cd}_1(q,u)$ is the function defined as
$$
\textrm{cd}_1=\textrm{cd}_1(q,u)=\frac{2\pi}{Kk_r}\sum^{\infty}_{n=0}\frac{(-1)^nq^{n+1/2}\cos\left((2n+3)\frac{\pi}{2K}u\right)}{1-q^{2n+1}}.\eqno{(98.1)}
$$
Hence we get:\\
\\
\textbf{Theorem 21.}\\
If we consider the nome 
\begin{equation}
q^a=iq^{1/2}e^{i\pi \theta_0/(2K)},
\end{equation}
which is parametrized as in Theorem 18, then
\begin{equation}
\textrm{cd}_1(q,\theta_0)=e^{-i\pi\theta_0/K}\textrm{cd}(q,\theta_0)-\frac{\sin\left(\frac{\pi\theta_0}{K}\right)}{m(q)K'}\left(\frac{d\theta(u)}{du}\right)_{u=a}
\end{equation}
and
\begin{equation}
\textrm{ss}(-q,m'(q)\theta_0)=\frac{1}{m(q) K'}\left(\frac{d\theta(u)}{du}\right)_{u=a}+i\cdot\textrm{cd}(q,\theta_0).
\end{equation}
\\

We now define $q=e(z)$, $z=x+iy$, $x,y$ reals with $y>0$, $-\frac{1}{2}<x<\frac{1}{2}$, $0<a<1$ and $\theta^{*}_0=\theta^{*}_{1}+i\theta^{*}_{2}$ with $\theta^{*}_{1}=(a-1+(2a-1)x)2K^{*}$, $\theta^{*}_{2}=(2a-1)y2K^{*}$. Then we have the next nome  
\begin{equation}
e^{\pi i a}q^a=-q^{1/2}e^{\pi i \theta^{*}_0/(2K^{*})}
\end{equation}
and holds the next relation (set $q\rightarrow -q$ in (100)):
$$
\textrm{cd}_1(-q,\theta^{*}_0)=e^{-i\pi\theta^{*}_0/(Km'(q))}\textrm{cn}\left(q,\frac{\theta^{*}_0}{m'(q)}\right)
-
$$
\begin{equation}
-\frac{2\pi i}{Km(q)}\sin\left(\frac{\pi\theta^{*}_0}{Km'(q)}\right)\sum^{\infty}_{n=0}\frac{q^{a(2n+1)}e^{\pi i a(2n+1)}}{1+q^{2n+1}},
\end{equation}
since
\begin{equation}
\textrm{cd}(-q,u)=\textrm{cn}\left(q,\frac{u}{m'(q)}\right).
\end{equation}

But also is
\begin{equation}
\textrm{ss}(q,u)=\textrm{cd}\left(-q, m'(q) u\right)\cot\left(\frac{\pi u}{K}\right)-\textrm{cd}_1\left(-q, m'(q) u\right)\csc\left(\frac{\pi u}{K}\right).
\end{equation} 
Hence from (105) and (103) we get
$$
\textrm{ss}\left(q,\frac{\theta_0^{*}}{m'(q)}\right)=
$$
$$
=\textrm{cn}\left(q,\frac{\theta_0^{*}}{m'(q)}\right)\cot\left(\frac{\pi \theta_0^{*}}{Km'(q)}\right)-\textrm{cd}_1\left(-q,\theta_0^{*}\right)\csc\left(\frac{\pi \theta_0^{*}}{Km'(q)}\right)=
$$
$$
=i\cdot\textrm{cn}\left(q,\frac{\theta_0^{*}}{m'(q)}\right)+\frac{2\pi i}{Km(q)}\sum^{\infty}_{n=0}\frac{q^{a(2n+1)}e^{\pi i a(2n+1)}}{1+q^{2n+1}}
$$
and we lead to the next theorems:\\
\\
\textbf{Theorem 22.}\\
If $\theta^{*}_0$ is of nome (102), then 
\begin{equation}
\textrm{ss}\left(q,\frac{\theta_0^{*}}{m'(q)}\right)=i\cdot\textrm{cn}\left(q,\frac{\theta_0^{*}}{m'(q)}\right)+\frac{2\pi i}{m(q)K}\sum^{\infty}_{n=0}\frac{q^{a(2n+1)}e^{\pi i a(2n+1)}}{1+q^{2n+1}}.
\end{equation}
\\
\textbf{Theorem 23.}\\
If we have the nome (99), then 
\begin{equation}
\frac{d\theta^{*}(a)}{da}=2\pi i (2z-1)\sum^{\infty}_{n=0}\frac{(-q)^{a(2n+1)}}{1+q^{2n+1}}.
\end{equation}
See also relations (120),(121) in Notes (ii) in Section 5 below.\\  
\\
\textbf{Theorem 24.}\\
Using the two nomes (99),(102) we get
\begin{equation}
\frac{\theta_0^{*}}{K^{*}}-\frac{\theta_0}{K}=2a-1,
\end{equation}
\begin{equation}
\frac{d\theta^{*}_0}{da}=2K^{*}+4zK^{*}\textrm{, }\frac{1}{m'(q)}\frac{d\theta_0^{*}}{da}=2K+2i K'.
\end{equation}
In the case $x=0$, $0<a<1$, there holds
\begin{equation}
\sum^{\infty}_{n=0}\frac{q^{a(2n+1)}e^{\pi i a(2n+1)}}{1+q^{2n+1}}=\sum^{\infty}_{n=0}\frac{(-q)^{a(2n+1)}}{1+q^{2n+1}}
\end{equation}
and
$$
\textrm{ss}\left(q,\frac{\theta_0^{*}}{m'(q)}\right)=i\cdot\textrm{cn}\left(q,\frac{\theta_0^{*}}{m'(q)}\right)+\frac{1}{m(-q)K(m'(-q))}\frac{d\theta^{*}(a)}{da}=
$$
\begin{equation}
=i\cdot\textrm{cn}\left(q,\frac{\theta_0^{*}}{m'(q)}\right)+\frac{1}{(1+2iy)m(q)K}\frac{d\theta^{*}(a)}{da}.
\end{equation}
See relations (118),(119) in Notes (i) of Section 5 below.

\section{Refinements and Applications}

\textbf{Application 1.}\\
If $y>0$ and $a=\frac{1}{2}$, $x=0$, $q=e(z)=e^{-2 \pi y}$ and $\theta_1=-K$, $\theta_2=0$. Then $\theta_0=\theta_1+i\theta_2=-K$, $K=K(m(q))$,   $A=iq^{1/2}e^{i\pi\theta_0/(2K)}=q^{1/2}$ and we have
\begin{equation}
\sum^{\infty}_{n=0}\frac{q^{n+1/2}}{1-q^{2n+1}}=\frac{-1}{4\pi y}\left(\frac{d\theta(a)}{da}\right)_{a=1/2}.
\end{equation} 
Hence from Corollary 2 with $y=\frac{\sqrt{r}}{2}$:
$$
\lim_{P\rightarrow K}\frac{\textrm{cd}_1\left(q,P\right)}{P-K}=1+\frac{2\pi^2}{k_rK^2}\sum_{n\geq0\textrm{, }n-odd}\frac{q^{n/2}}{1-q^n}=
$$
\begin{equation}
=1-\frac{\pi}{2m(q)yK^2}\left(\frac{d\theta(a)}{da}\right)_{a=1/2}.
\end{equation} 
\\
\textbf{Application 2.}\\
If $a=1$, $x=0$, $y>0$, then $q=e^{-2\pi y}$ and
\begin{equation}
\sum^{\infty}_{n=0}\frac{1}{e^{2 (2n+1)\pi y}-1}=\frac{-1}{4 \pi y}\left(\frac{d\theta(a)}{da}\right)_{a=1}.
\end{equation} 
\\
\textbf{Application 3.}\\
Set $a=\frac{1}{2}$, $x=0$, $y>0$, then $q=e^{-2\pi y}$ and $\theta^{*}_0=-K^{*}=-K(m(-q))$. Hence from Theorem 24 relation (111) we have
\begin{equation}
\textrm{ss}\left(q,\frac{-K^{*}}{m'(q)}\right)=i\cdot \textrm{cn}\left(q,\frac{K^{*}}{m'(q)}\right)+\frac{1}{m(-q)K(m'(-q))}\left(\frac{d\theta^{*}(a)}{da}\right)_{a=1/2}
\end{equation}
But it holds $K^{*}=K(m(q))m'(q)$,  $(K')^{*}=K(m'(-q))=-i(1+2z)K(m(-q))=-i(1+2iy)K(m(q))m'(q)$ and $m(-q)=i\frac{m(q)}{m'(q)}$. Hence it must be $m(-q)K^{*}{'}=-i(1+2iy)K(m(q))m'(q)i\frac{m(q)}{m'(q)}=(1+2iy)K(m(q))m(q)$. Consequently
$$
-\textrm{ss}\left(q,K\right)=i\cdot \textrm{cn}\left(q,K\right)+\frac{1}{(1+2iy) m(q)K}\left(\frac{d\theta^{*}(a)}{da}\right)_{a=1/2}.
$$
Or equivalently from $\textrm{cn}(q,K)=0$:
\begin{equation}
\textrm{ss}(q,K)=\frac{-1}{(1+2iy)m(q)K}\left(\frac{d\theta^{*}(a)}{da}\right)_{a=1/2}
\end{equation}
and
\begin{equation}
\sum^{\infty}_{n=0}\frac{(-1)^nq^{n+1/2}}{1+q^{2n+1}}=\frac{-1}{2\pi(1+2iy)}\left(\frac{d\theta^{*}(a)}{da}\right)_{a=1/2}\textrm{, }q=e^{-2\pi y}\textrm{, }y>0.
\end{equation}
\\
\textbf{Notes.}\\
\textbf{i)} In case of $x=0$, $0<a<1$, we can prove relation (111) using the arguments of Application 3. Hence then, $\theta^{*}_0=-2\left((1-a)+i(1-2a)y\right)K^{*}$ and equivalently 
$$
\widetilde{h}^{*}_0:=\frac{\theta^{*}_0}{m'(q)}=-2\left[(1-a)+i(1-2a)y\right]K=
$$
\begin{equation}
=2(a-1)K+i(2a-1)K',
\end{equation} 
with
\begin{equation}
\textrm{ss}\left(q,\widetilde{h}^{*}_0\right)=i\cdot \textrm{cn}\left(q,\widetilde{h}^{*}_0\right)+\frac{1}{(1+2iy)m(q) K}\frac{d\theta^{*}(a)}{da},
\end{equation}
where $q=e^{-2\pi y}$, $y>0$.\\
\textbf{ii)} When $q=e(z)$, $z=x+iy$, $y>0$, $-\frac{1}{2}<x+\frac{1}{2}<\frac{1}{2}$, $0<a<1$, we have (general case)
$$
h_0^{*}=\frac{\theta_0^{*}}{m'(q)}=2((a-1)+(2a-1)x)K+i(2a-1)yK'=
$$
\begin{equation}
=2(a-1)K+i(2a-1)K'
\end{equation}
and
$$
\frac{d\theta^{*}(a)}{da}= 2i\pi\frac{K\left(m'(-q)\right)m(-q)}{K\left(m(q)\right)m(q)}\sum^{\infty}_{n=0}\frac{(-q)^{a(2n+1)}}{1+q^{2n+1}}=
$$
\begin{equation}
=2\pi i (1+2z)\sum^{\infty}_{n=0}\frac{(-q)^{a(2n+1)}}{1+q^{2n+1}}.
\end{equation}
\textbf{iii)} 
Taking the logarithms in both sides of (81) and expanding them, then differentiating with respect to $x$, we get after using divisor sums
\begin{equation}
\left(\frac{\partial\theta(q,u)}{\partial u}\right)_{u=a}=\frac{d\theta(a)}{da}=\frac{2q^a\log(q)}{1-q^{2a}}+2\log(q)\sum^{\infty}_{n=1}q^n\sum_{\scriptsize{
\begin{array}{cc}
	d|n\\
	d>a\\
	n/d-odd 
\end{array}
}}1.
\end{equation}
This relation holds for $a$ positive integer. Hence from (127) below we get ($q=e^{-2\pi y}$, $y>0$):
\begin{equation}
\sum^{\infty}_{n=0}\frac{1}{e^{2(2n+1)\pi y}-1}=\frac{q}{1-q^{2}}+\sum^{\infty}_{n=1}q^n\sum_{\scriptsize{
\begin{array}{cc}
	d|n\\
	d>1\\
	n/d-odd 
\end{array}
}}1.
\end{equation}
\\
\textbf{Application 4.}\\
Let $a$ be positive integer, then
\begin{equation}
\left(\frac{\partial\theta(q,u)}{\partial u}\right)_{u=a+1}-\left(\frac{\partial\theta(q,u)}{\partial u}\right)_{u=a}=-\frac{2q^a}{1-q^{2a}}\log(q).
\end{equation}
Hence easily with $a$ positive integer:
\begin{equation}
\left(\frac{\partial\theta(q,u)}{\partial u}\right)_{u=a}=-2\log(q)\sum^{a-1}_{n=1}\frac{q^n}{1-q^{2n}}+2\log(q)\sum^{\infty}_{n=1}\frac{q^n}{1-q^{2n}}
\end{equation}
and if $q=e^{-2\pi y}$, $y>0$, then
\begin{equation}
\frac{1}{2\pi y}\frac{d\theta\left(a\right)}{da}=\sum^{a-1}_{n=1}\frac{1}{\sinh(2\pi n y)}-\sum^{\infty}_{n=1}\frac{1}{e^{2\pi ny}-1}-\sum^{\infty}_{n=1}\frac{1}{e^{2\pi n y}+1}.
\end{equation}
Hence also
\begin{equation}
\frac{1}{2\pi y}\frac{d\theta\left(a\right)}{da}=\sum^{a-1}_{n=1}\frac{1}{\sinh(2\pi n y)}-2\sum^{\infty}_{n=0}\frac{1}{e^{2\pi (2n+1)y}-1}.
\end{equation}
\\
\textbf{Application 5.}\\
If the nome is $q=e(z)$, $Im(z)>0$ and $a=\frac{1}{2}$, then $\theta_0=-K$, and from (87), for all $|q|<1$ we have:
\begin{equation}
\theta\left(q,\frac{1}{2}\right)=-\frac{1}{2}\log\left(\frac{m'(q)}{1+m(q)}\right).
\end{equation}
\\
\textbf{Application 6.}\\
An example of evaluation using the modular angle  $\theta(q,a)$ function is: Setting $\lambda=\lambda_1+i\lambda_2$ in Theorem 10, then using Theorem 21 with $x=0$, $y=\frac{\sqrt{r}}{2}$, $a=\frac{1+\lambda_1}{2}$, $r=\frac{1}{\lambda_2^2}$, $z=x+i y$, $q=e(z)$, $\theta_0=\theta_1+i\theta_2$, $\theta_1=-K$, $\theta_2=(2a-1)\sqrt{r}K$, after simplifications, we get the next formula:
$$
\sum^{\infty}_{n=0}\frac{(-1)^ne^{-\pi(n+1/2)\lambda/\lambda_2}}{\sinh\left((n+1/2)\pi/\lambda_2\right)}=\frac{i\lambda_2\theta^{(0,1)}\left(q,\frac{1+\lambda_1}{2}\right)}{\pi\sqrt{r}}.\eqno{(128.1)}
$$
Where $|\lambda|<1$, $k=m=m\left(e^{-\pi/\lambda_2}\right)$, $K=K(k)$, $K'=K(k')$, $k'=\sqrt{1-k^2}$ and $\theta^{(0,1)}(q,a)=\frac{d\theta(q,a)}{da}$.\\
Hence we get evaluations of the kind: If $\lambda_1=2a-1$, $\lambda_2=1/\sqrt{r}$, $r>0$, then:
\begin{equation}
\sum^{\infty}_{n=0}\frac{ e^{-\pi(n+1/2)(2a-1)\sqrt{r}}}{\sinh\left((n+1/2)\pi\sqrt{r}\right)}=\frac{-\theta^{(0,1)}\left(q,a\right)}{\pi\sqrt{r}}\textrm{, where }(2a-1)^2+1/r<1.
\end{equation}
For $a$ positive integer and $r>0$, we have (in view of (127)): 
$$
\sum^{\infty}_{n=0}\frac{ e^{-\pi(n+1/2)(2a-1)\sqrt{r}}}{\sinh\left((n+1/2)\pi\sqrt{r}\right)}=-\sum^{a-1}_{n=1}\frac{1}{\sinh(\pi n\sqrt{r})}+2\sum^{\infty}_{n=0}\frac{1}{e^{\pi(2n+1)\sqrt{r}}-1}.\eqno{(129.1)}
$$ 
\\
\textbf{Application 7.}\\
If $q=e^{-2\pi y}$, $y>0$, then
$$
\sum^{\infty}_{n=1}\frac{q^{n/2}}{1-q^n}=\sum^{\infty}_{n=0}\frac{q^{n+1/2}}{1-q^{2n+1}}+\sum^{\infty}_{n=0}\frac{1}{e^{2\pi (2n+1)y}-1}.
$$ 
Hence from (114),(112) we get
\begin{equation}
-4\pi y\sum^{\infty}_{n=1}\frac{q^{n/2}}{1-q^n}=\left(\frac{d\theta(a)}{da}\right)_{a=1/2}+\left(\frac{d\theta(a)}{da}\right)_{a=1}.
\end{equation}

\section{Appendix: Evaluations of $U$ using theta functions and series with divisor sums coefficients}

From [12] we have\\
\\
\textbf{Theorem A1.}\\
If $|q|<1$, $a,p$ real and $p>0$, then
\begin{equation}
-1+\frac{2}{1-U\left(q^{k+h},-q^{k-h},q^{2k}\right)}=\frac{\sum^{\infty}_{n=-\infty}q^{kn^2+hn}}{\sum^{\infty}_{n=-\infty}(-1)^nq^{kn^2+hn}},
\end{equation}
or equivalently the next continued fraction expansion holds
\begin{equation}
\frac{\sum^{\infty}_{n=-\infty}q^{pn^2/2+(p-2a)n/2}}{\sum^{\infty}_{n=-\infty}(-1)^nq^{pn^2/2+(p-2a)n/2}}=-1+\frac{2}{1-U\left(q^{a},-q^{p-a},q^{p}\right)}.
\end{equation}
\\

We now define (see [26]):
\begin{equation}
\left[a,p;q\right]^{-}_{\infty}:=\prod^{\infty}_{n=0}\left(1-q^{p n+a}\right)\left(1-q^{p n+p-a}\right)
\end{equation}
and
\begin{equation}
\left[a,p;q\right]^{+}_{\infty}:=\prod^{\infty}_{n=0}\left(1+q^{p n+a}\right)\left(1+q^{p n+p-a}\right)
\end{equation}
We call $[a,p;q]^{\pm}_{\infty}$ agiles. Multiplying the two ''agiles'' we have
\begin{equation}
\left[a,p;q\right]^{-}_{\infty}\left[a,p;q\right]^{+}_{\infty}=\left[2a,2p;q\right]^{-}_{\infty}=\left[a,p;q^2\right]^{-}_{\infty}.
\end{equation}
But in view of (64),(66) and
\begin{equation}
\left(-1+\frac{2}{1-U(q^a,-q^b;q^c)}\right)^2=\left(-1+\frac{2}{1-u_0(q^a;q^c)}\right)\left(-1+\frac{2}{1-u_0(q^b;q^c)}\right),
\end{equation}
we have
$$
\left(\frac{\left[a,p;q\right]^{-}_{\infty}}{\left[a,p;q\right]^{+}_{\infty}}\right)^2=\left(\frac{\left(q^a;q^p\right)_{\infty}}{\left(-q^a;q^p\right)_{\infty}}\right)^2\left(\frac{\left(q^{p-a};q^p\right)_{\infty}}{\left(-q^{p-a};q^p\right)_{\infty}}\right)^2=
$$
\begin{equation}
=\left(-1+\frac{2}{1-U\left(q^a,-q^{p-a};q^p\right)}\right)^{-2}=\left(\frac{\sum^{\infty}_{n=-\infty}q^{pn^2/2+(p-2a)n/2}}{\sum^{\infty}_{n=-\infty}(-1)^nq^{pn^2/2+(p-2a)n/2}}\right)^{-2},
\end{equation}
where
\begin{equation}
\vartheta_3(a,b;q):=\sum^{\infty}_{n=-\infty}q^{an^2+bn}\textrm{, }|q|<1
\end{equation}
and
\begin{equation}
\vartheta_4(a,b;q):=\sum^{\infty}_{n=-\infty}(-1)^nq^{an^2+bn}\textrm{, }|q|<1.
\end{equation}
Hence from Theorem A1 (relation (131)) we get
\begin{equation}
\frac{\left[a,p;q\right]^{-}_{\infty}}{\left[a,p;q\right]^{+}_{\infty}}=\frac{\vartheta_4\left(\frac{p}{2},\frac{p-2a}{2};q\right)}{\vartheta_3\left(\frac{p}{2},\frac{p-2a}{2};q\right)}.
\end{equation}
Relations (135) and (140) lead us to the next\\ 
\\
\textbf{Theorem A2.}\\
If $|q|<1$, then
\begin{equation}
\left[a,p;q^2\right]^{-}_{\infty}=\frac{\vartheta_3\left(\frac{p}{2},\frac{p-2a}{2};q\right)}{\vartheta_4\left(\frac{p}{2},\frac{p-2a}{2};q\right)}\left(\left[a,p,q\right]^{-}_{\infty}\right)^2
\end{equation}
and
\begin{equation}
\left[a,p;q^2\right]^{-}_{\infty}=\frac{\vartheta_4\left(\frac{p}{2},\frac{p-2a}{2};q\right)}{\vartheta_3\left(\frac{p}{2},\frac{p-2a}{2};q\right)}\left(\left[a,p,q\right]^{+}_{\infty}\right)^2
\end{equation}
\\
\textbf{Example.}\\
Set 
\begin{equation}
G(q):=\sum^{\infty}_{n=0}\frac{q^{n^2}}{(q;q)_n}=\frac{1}{[1,5;q]^{-}_{\infty}},
\end{equation}
then
\begin{equation}
G(q^2)=\frac{\vartheta_4\left(\frac{5}{2},\frac{-3}{2};q\right)}{\vartheta_3\left(\frac{5}{2},\frac{-3}{2};q\right)}G(q)^2=\frac{\vartheta_4\left(\frac{5}{2},\frac{3}{2};q\right)}{\vartheta_3\left(\frac{5}{2},\frac{3}{2};q\right)}G(q)^2
\end{equation}
Also set
\begin{equation}
H(q):=\sum^{\infty}_{n=1}\frac{q^{n^2+n}}{(q;q)_n}=\frac{1}{[2,5;q]^{-}_{\infty}},
\end{equation}
then
\begin{equation}
H(q^2)=\frac{1}{\left[2,5;q^2\right]^{-}_{\infty}}=\frac{\vartheta_4\left(\frac{5}{2},\frac{-1}{2};q\right)}{\vartheta_3\left(\frac{5}{2},\frac{-1}{2};q\right)}H(q)^2=\frac{\vartheta_4\left(\frac{5}{2},\frac{1}{2};q\right)}{\vartheta_3\left(\frac{5}{2},\frac{1}{2};q\right)}H(q)^2.
\end{equation}
But the Rogers-Ramanujan continued fraction is
\begin{equation}
R(q)=q^{1/5}\frac{H(q)}{G(q)}=q^{1/5}\frac{\left[1,5;q\right]^{-}_{\infty}}{\left[2,5;q\right]^{-}_{\infty}}.
\end{equation}
Hence
$$
R(q^2)=q^{2/5}\frac{\vartheta_4\left(\frac{5}{2},\frac{1}{2};q\right)}{\vartheta_3\left(\frac{5}{2},\frac{1}{2};q\right)}\frac{\vartheta_3\left(\frac{5}{2},\frac{3}{2};q\right)}{\vartheta_4\left(\frac{5}{2},\frac{3}{2};q\right)}\left(\frac{H(q)}{G(q)}\right)^2=
$$
\begin{equation}
=\frac{R(q)^2}{q^{-1/5}R(q)}\frac{\vartheta_3\left(\frac{5}{2},\frac{3}{2};q\right)}{\vartheta_3\left(\frac{5}2{},\frac{1}{2};q\right)}=q^{1/5}R(q)\frac{\vartheta_3\left(\frac{5}{2},\frac{3}{2};q\right)}{\vartheta_3\left(\frac{5}{2},\frac{1}{2};q\right)}.
\end{equation}
Hence we get the next evaluation
\begin{equation}
\frac{\vartheta_3\left(\frac{5}{2},\frac{3}{2};q\right)}{\vartheta_3\left(\frac{5}{2},\frac{1}{2};q\right)}=q^{-1/5}\frac{R(q^2)}{R(q)}\textrm{, }|q|<1.
\end{equation}

Generalizing the above arguments for Ramanujan quantities we can write\\
\\
\textbf{Theorem A3.}\\
If $a,b$ positive odd integers and $p$ positive even integer such $a<b<a+b<p$, then
\begin{equation}
R(a,b,p;-q)=\frac{\vartheta_3\left(\frac{p}{2},\frac{p-2a}{2};q\right)}{\vartheta_3\left(\frac{p}{2},\frac{p-2b}{2};q\right)}.
\end{equation} 
\\
\textbf{Notes.}\\
i) The above theorem is based on the next argument, which is elementary to prove:\\
If $a,b$ are odd positive integers and $p$ in even positive integer and $a<b<a+b<p$, then
\begin{equation}
R(a,b,p;q^2)=R(a,b,p;q)R(a,b,p;-q).
\end{equation}
ii) On the other hand it have been proved in [26] that
\begin{equation}
R\left(a,b,p;q\right)=\frac{\vartheta_4\left(\frac{p}{2},\frac{p-2a}{2};q\right)}{\vartheta_4\left(\frac{p}{2},\frac{p-2b}{2};q\right)}\textrm{, }|q|<1,
\end{equation}
where
\begin{equation}
R(a,b,p;q)=\frac{\left[a,p;q\right]^{-}_{\infty}}{\left[b,p;q\right]^{-}_{\infty}}.
\end{equation}
\\
\textbf{Theorem A4.}(see [27])\\
If $|q|<1$ and $a>0$, then 
\begin{equation}
\left[a,p,q\right]^{+}_{\infty}=\frac{\vartheta_3\left(\frac{p}{2},\frac{p-2a}{2};q\right)}{f\left(-q^p\right)}
\end{equation}
and
\begin{equation}
\left[a,p;q\right]^{-}_{\infty}=\frac{\vartheta_4\left(\frac{p}{2},\frac{p-2a}{2};q\right)}{f\left(-q^p\right)},
\end{equation}
where
\begin{equation}
f(-q)=\prod^{\infty}_{n=1}\left(1-q^n\right).
\end{equation}
\\

Using Theorems A1,A4, we get
$$
\left[a,p;q\right]^{-}_{\infty}=\left(-1+\frac{2}{1-U\left(q^a,-q^{p-a};q^p\right)}\right)^{-2}\left[a,p;q\right]^{+}_{\infty}.
$$
Hence\\
\\
\textbf{Theorem A5.}\\
If $0<a<p$ and $|q|<1$, then holds the following continued fraction expansion
\begin{equation}
\frac{\vartheta_4\left(\frac{p}{2},\frac{p}{2}-a;q\right)}{\vartheta_3\left(\frac{p}{2},\frac{p}{2}-a;q\right)}=\left(-1+\frac{2}{1-U\left(q^a,-q^{p-a};q^p\right)}\right)^{-2}
\end{equation}
Moreover if $a,p$ are positive integers, then
\begin{equation}
\log\left(-1+\frac{2}{1-U\left(q^a,-q^{p-a};q^p\right)}\right)=4\sum^{\infty}_{n=1}q^n\sum_{\scriptsize
\begin{array}{cc}
AB=n\\
A\equiv1(2)\\
B\equiv \pm a(p)	
\end{array}}\frac{1}{A}.
\end{equation}
\\
\textbf{Proof.}\\
Using relation (67) we get 
$$
-2\log\left(-1+\frac{2}{1-U(q^a,-q^{p-a};q^p)}\right)=
-2\log\left(-1+\frac{2}{1-u_0\left(q^a,q^p\right)}\right)-
$$
$$
-2\log\left(-1+\frac{2}{1-u_0\left(q^{p-a},q^p\right)}\right)=
-8\sum^{\infty}_{n\geq 0,n-odd}\frac{q^{na}}{n(1-q^{pn})}-
$$
$$
-8\sum^{\infty}_{n\geq 0,n-odd}\frac{q^{n(p-a)}}{n(1-q^{pn})}
$$
But it holds
$$
-8\sum^{\infty}_{n\geq 0,n-odd}\frac{q^{na}}{n(1-q^{np})}-8\sum^{\infty}_{n\geq 0,n-odd}\frac{q^{n(p-a)}}{n(1-q^{pn})}
=-8\sum_{n,l\geq0,n-odd}\frac{q^{na+lnp}}{n}-
$$
$$
-8\sum_{n,l\geq0,n-odd}\frac{q^{n(p-a)+lnp}}{n}=-8\sum^{\infty}_{n=1}q^n\sum_{\scriptsize
\begin{array}{cc}
AB=n\\
A\equiv1(2)\\
B\equiv a(p)	
\end{array}}\frac{1}{A}-8\sum^{\infty}_{n=1}q^n\sum_{\scriptsize
\begin{array}{cc}
AB=n\\
A\equiv 1(2)\\
B\equiv -a(p)	
\end{array}}\frac{1}{A}. 
$$ 
$qed$.\\

In the same way we prove\\
\\
\textbf{Theorem A6.}\\
Let $a,b,p$ be positive integers with $a,b<p$ and $|q|<1$. If $\epsilon=\pm 1$, then
$$
\log\left(-1+\frac{2}{1-U\left(q^a,\epsilon  q^b;q^p\right)}\right)
=2\sum^{\infty}_{n=1}q^n\sum_{\scriptsize
\begin{array}{cc}
AB=n\\
A\equiv1(2)\\
B\equiv a(p)	
\end{array}}\frac{1}{A}-
$$
\begin{equation}
-2\epsilon\sum^{\infty}_{n=1}q^n\sum_{\scriptsize
\begin{array}{cc}
AB=n\\
A\equiv1(2)\\
B\equiv b(p)	
\end{array}}\frac{1}{A}.
\end{equation}
\\
\textbf{Corollary.}\\
If $a,p$ are positive integers with $a<p$ and $|q|<1$, then
\begin{equation}
\log\vartheta_3\left(\frac{p}{2},\frac{p-2a}{2};q\right)-\log\vartheta_4\left(\frac{p}{2},\frac{p-2a}{2};q\right)=2\sum^{\infty}_{n=1}q^n\sum_{\scriptsize
\begin{array}{cc}
AB=n\\
A\equiv1(2)\\
B\equiv \pm a(p)	
\end{array}}\frac{1}{A}. 
\end{equation}
\\
\textbf{Notes.}\\
If $a,p$ positive integers and $|q|<1$, then\\
\textbf{i)}
\begin{equation}
\log\vartheta_3\left(\frac{p}{2},\frac{p}{2}-a;q\right)=\sum^{\infty}_{n=1}\frac{q^{pn}}{n}\sigma_1(n)+\sum^{\infty}_{n=1}q^n\sum_{\scriptsize
\begin{array}{cc}
AB=n\\
B\equiv \pm a(p)	
\end{array}}\frac{(-1)^A}{A},
\end{equation}
where $\sigma_1(n)$ is the sum of positive divisors of $n$.\\
\textbf{ii)}
\begin{equation}
\log\left(\left[a,p;q\right]^{-}_{\infty}\right)=-\sum^{\infty}_{n=1}q^n\sum_{\scriptsize
\begin{array}{cc}
AB=n\\
B\equiv \pm a(p)	
\end{array}}\frac{1}{A}
\end{equation}
and
\begin{equation}
\log\left(\left[a,p;q\right]^{+}_{\infty}\right)=-\sum^{\infty}_{n=1}q^n\sum_{\scriptsize
\begin{array}{cc}
AB=n\\
B\equiv \pm a(p)	
\end{array}}\frac{(-1)^A}{A}.
\end{equation}
\textbf{iii)}
\begin{equation}
\log\left(-1+\frac{2}{1-U(q^a,-q^b,q^p)}\right)=4\sum^{\infty}_{n=1}q^n\sum_{\scriptsize
\begin{array}{cc}
AB=n\\
A\equiv1(2)\\
B\equiv a,b(p)	
\end{array}}\frac{1}{A}
\end{equation}
\textbf{iv)} 
\begin{equation}
\log\left(\left(\pm q^a;q^p\right)_{\infty}\right)=-\sum^{\infty}_{n=1}q^n\sum_{\scriptsize
\begin{array}{cc}
AB=n\\
B\equiv a(p)	
\end{array}}\frac{(\pm 1)^A}{A}
\end{equation}
\\
\textbf{Theorem A7.}\\
Let $a,b,p$ be positive integers such $a,b<p$ and $-1<q<1$, $\epsilon=\pm 1$, then   
$$
\log\left(-1+\frac{2}{1-U\left(q^a,\epsilon q^b;q^p\right)}\right)+\log\left(-1+\frac{2}{1-U\left(q^{p-a},\epsilon q^{p-b};q^p\right)}\right)=
$$
\begin{equation}
=\log\left(\frac{\vartheta_3\left(\frac{p}{2},\frac{p-2a}{2};q\right)}{\vartheta_4\left(\frac{p}{2},\frac{p-2a}{2};q\right)}\right)-\epsilon \log\left(\frac{\vartheta_3\left(\frac{p}{2},\frac{p-2b}{2};q\right)}{\vartheta_4\left(\frac{p}{2},\frac{p-2b}{2};q\right)}\right).
\end{equation}
\\
\textbf{Proof.}\\
Straight forward evaluation using Theorem A6.\\
\\
\textbf{Theorem A8.}\\
If $|q|<1$, $x,y$ real and $a$, $b$, $p$ integers such $0<a,b<p$, then
\begin{equation}
\log\left(-1+\frac{2}{1-U\left(xq^a,-xq^b,q^p\right)}\right)=\sum^{\infty}_{n=1}q^n\sum_{\scriptsize
\begin{array}{cc}
	AB=n\\
	A\equiv 1(2)\\
	B\equiv a,b(p)
\end{array}
\normalsize}\frac{x^A}{A}
\end{equation}
and
$$
\log\left(-1+\frac{2}{1-U\left(xq^a,-yq^b,q^p\right)}\right)=\sum^{\infty}_{n=1}q^n\sum_{\scriptsize
\begin{array}{cc}
	AB=n\\
	A\equiv 1(2)\\
	B\equiv a(p)
\end{array}
\normalsize}\frac{x^A}{A}+
$$
\begin{equation}
+\sum^{\infty}_{n=1}q^n\sum_{\scriptsize
\begin{array}{cc}
	AB=n\\
	A\equiv 1(2)\\
	B\equiv b(p)
\end{array}
\normalsize}\frac{y^A}{A}.
\end{equation}
\\
\textbf{Proof.}\\
Use the identity
\begin{equation}
\log\left(\frac{\left(-xq^a;q^p\right)_{\infty}}{\left(xq^a;q^p\right)_{\infty}}\right)=2\sum^{\infty}_{n=1}q^n\sum_{\scriptsize
\begin{array}{cc}
	AB=n\\
	A\equiv 1(2)\\
	B\equiv a(p)
\end{array}
\normalsize}\frac{x^A}{A},
\end{equation}
along with (64),(65),(66) and
\begin{equation}
\left(-1+\frac{2}{1-U(a,-b;q)}\right)^2=\left(-1+\frac{2}{1-u_0(a,q)}\right)\left(-1+\frac{2}{1-u_0(b,q)}\right).
\end{equation}
$qed$\\
\\
Having proved the above theorem we can easy state the following main theorem for the evaluation of $U$ continued fraction:\\
\\
\textbf{Main Theorem A1.}\\
If $|q|<1$ and $x,y\in\textbf{R}$ such $|xq|<1$, $|yq|<1$, then
$$
\log\left(-1+\frac{2}{1-U\left(xq,yq,q\right)}\right)=2\sum^{\infty}_{n=1}q^n\sum_{\scriptsize
\begin{array}{cc}
	d|n\\
	d-odd\\
\end{array}
\normalsize}\frac{x^d}{d}-
$$
\begin{equation}
-2\sum^{\infty}_{n=1}q^n\sum_{\scriptsize
\begin{array}{cc}
	d|n\\
	d-odd\\
\end{array}
\normalsize}\frac{y^d}{d}=2\sum_{\scriptsize
\begin{array}{cc}
	n\geq 1\\
	n-odd\\
\end{array}
\normalsize}\frac{x^nq^n}{n(1-q^n)}-2\sum_{\scriptsize
\begin{array}{cc}
	n\geq 1\\
	n-odd\\
\end{array}
\normalsize}\frac{y^nq^n}{n(1-q^n)}.
\end{equation}

\newpage

\centerline{\bf References}\vskip .2in

\noindent

[1]: M.Abramowitz and I.A.Stegun, 'Handbook of Mathematical Functions'. Dover Publications, New York., (1972).

[2]: B.C. Berndt, 'Ramanujan`s Notebooks Part II'. Springer Verlang, New York., (1989).

[3]: B.C. Berndt, 'Ramanujan`s Notebooks Part III'. Springer Verlang, New York., (1991).

[4]: I.S. Gradshteyn and I.M. Ryzhik, 'Table of Integrals, Series and Products'. Academic Press., (1980).

[5]: L. Lorentzen and H. Waadeland, Continued Fractions with Applications. Elsevier Science Publishers B.V., North Holland., (1992).  

[6]: H.S. Wall. 'Analytic Theory of Continued Fractions'. Chelsea Publishing Company, Bronx, N.Y., (1948). 

[7]: E.T. Whittaker and G.N. Watson. 'A course on Modern Analysis'. Cambridge U.P., (1927).

[8]: J.V. Armitage, W.F. Eberlein. 'Elliptic Functions'. Cambridge University Press., (2006).

[9]: J.M. Borwein, M.L. Glasser, R.C. McPhedran, J.G. Wan, I.J. Zucker. 'Lattice Sums Then and Now'. Cambridge University Press. New York., (2013).

[10]: J.M. Borwein and P.B. Borwein. 'Pi and the AGM: A Study in Analytic Number Theory and Computational Complexity', Wiley, New York., (1987).

[11]: M.L. Glasser and N.D. Bagis. 'Some Applications of the Poisson Summation Formula'. arXiv:0812.0990, (2008)  

[12]: N.D. Bagis and M.L. Glasser. 'Evaluations of a Continued Fraction of Ramanujan'. Rend. Sem. Mat. Univ. Padova. Vol 133., (2015). 

[13]: S.C. Milne. 'Infinite Families of Exact Sums of Squares Formulas, Jacobi Elliptic Functions, Continued Fractions, and Schur Functions'.\\ arXiv:math/0008068v2 [math.NT] 7 Juan. 2001.

[14]: N.D. Bagis. 'On certain theta functions and modular forms in Ramanujan theories'. arXiv:1511.03716v2 [math.GM] 6 Dec 2017. 

[15]: N.D. Bagis, M.L. Glasser. 'Conjectures on the Evaluation of Alternative Modular Bases and Formulas Approximating $1/\pi$'. Journal of Number Theory. (Elsevier), (2012).

[16]: N.D. Bagis, M.L. Glasser. 'Conjectures on the evaluation of certain functions with algebraic properties'. Journal of Number Theory. 155 (2015), 63-84

[17]: N.D. Bagis. 'On the Complete Evaluation of Jacobi Theta Functions'. arXiv:1503.01141v1 [math.GM],(2015).

[18]: Don. Zagier. 'Elliptic Modular Forms and Their Applications'. Available from zagier@mpim-bonn.mpg.de

[19]: Toshijune Miyake. 'Modular Forms'. Springer Verlang, Berlin, Heidelberg, (1989).

[20]: Carlos J. Moreno, Samuel S. Wagstaff. Jr. 'Sums of Squares of Integers'. Chapman and Hall/CRC, Taylor and Francis Group, (2006).

[21]: Shaun Cooper, Dongxi Ye. 'Level 14 and 15 analogues of Ramanujan's elliptic functions to alternative bases'. Transactions of the American Mathematical Society, \textbf{368} (2016), 7883-7910.

[22]: I.J.Zucker. 'The summation of series of hyperbolic functions'. SIAM J. Math. Ana.10.192, 1979.

[23]: D. Broadhurst. 'Solutions by Radicals at Singular Values $k_N$ from New Class Invariants for $N\equiv3mod8$'. arXiv:0807.2976 (math-phy), (2008) 

[24]: S. Yakubovich, P. Drygas and V. Mityushev. 'Closed-form evaluation of two-dimensional static lattice sums'. Proc. R. Soc. A 472:20160510, (2017)  

[25]: N. Bagis and M.L. Glasser. 'Ramanujan type $1/\pi$ approximation formulas'. Journal of Number Theory. 133 (2013) 3453-3469. Elsevier.
 
[26]: N.D. Bagis. 'Generalizations of Ramanujan's Continued Fractions'. arXiv:1107.2393v2[math.GM] 7 Aug 2012.

[27]: N.D. Bagis, M.L. Glasser. 'Jacobian Elliptic Functions, Continued Fraction and Ramanujan Quantities'. arXiv:1001.2660v1 [math.GM] 15 Jan 2010.

\end{document}